\documentclass[seceqn]{elsart3p}

\usepackage{amsfonts,amsmath,xy}
\usepackage{latexsym,lineno}

\journal{Systems \& Control Letters}
\newtheorem{theorem}{Theorem}[section]
\newtheorem{proposition}[theorem]{Proposition}

\newtheorem{lemma}[theorem]{Lemma}

\newtheorem{definition}{Definition}[section]

\newtheorem{remark}{Remark}[section]

\newcommand{\R}{\ensuremath{\mathbb{R}}}

\def \e {\varepsilon}

\def \l {\lambda}

\def \eqref {\ref}

\begin{document}
\begin{frontmatter}

\title{An efficient algorithm for positive realizations}

\author{Wojciech Czaja\thanksref{harp}\thanksref{eu}},
\thanks[harp]{Partially supported by {\it European Commission}
Harmonic Analysis and Related Problems 2002-2006 IHP Network,
Contract Number: HPRN-CT-2001-00273 - HARP.}
\thanks[eu]{Partially supported by {\it European Commission} grant MEIF-CT-2003-500685.}
\ead{wojtek@math.umd.edu}
\address{Institute of Mathematics, University of Wroclaw, Pl. Grunwaldzki 2/4, 50-384
Wroclaw, POLAND,\\ and Department of Mathematics, University of Maryland, College Park, MD 20742, USA}

\author{Philippe Jaming\thanksref{harp}\thanksref{balaton}},
\thanks[balaton]{Partially supported by the Hungarian-French Scientific and
Technological Governmental Cooperation, no. F-10/04}
\ead{philippe.jaming@univ-orleans.fr}
\address{MAPMO-F\'ed\'eration Denis Poisson, Universit\'e d'Orl\'eans,
B.P. 6759, 45067 Orl\'eans cedex 2, FRANCE}
\author{M\'at\'e Matolcsi\corauthref{cor}\thanksref{harp}\thanksref{balaton}\thanksref{hung}}
\corauth[cor]{Corresponding author.}
\thanks[hung]{Partially supported by OTKA-T047276, T049301, PF64601.}
\ead{matomate@renyi.hu}
\address{Alfr\'ed R\'enyi Institute of Mathematics, Budapest, H-1053, HUNGARY}

\begin{abstract}
We observe that successive applications of known results from the
theory of positive systems lead to an {\it efficient general
algorithm} for positive realizations of transfer functions. We
give two examples to illustrate the algorithm, one of which
complements an earlier result of \cite{large}. Finally, we improve
a lower-bound of \cite{mn2} to indicate that the algorithm is
indeed efficient in general.
\end{abstract}


\end{frontmatter}

\section{Introduction}
Given the transfer function
$$
H(z)=\frac{p_{1}z^{n-1}+...+p_n}{z^n+q_{1}z^{n-1}+...+q_n}, \ \ \
p_j, q_j\in\R ,
$$
of a {\sl discrete time-invariant linear
SISO system} of McMillan degree $n$, we say that a triple
$\mathbf{A}\in{\R}^{n\times n}$, $\mathbf{b},\mathbf{c}\in{\R}^n$
is an {\it $n$th order realization of $H(z)$} if it satisfies the
condition:
$$
H(z)=\mathbf{c}^T(z\mathbf{I}-\mathbf{A})^{-1}\mathbf{b}.
$$
It is known that an $n$th order realization of $H(z)$
always exists (see, e.g. \cite[Chapter 9]{farri}). In this note,
however, we are interested in the {\sl positive realization
problem}, i.e. finding $\mathbf{A}$, $\mathbf{b}$, $\mathbf{c}$
with nonnegative entries (and possibly of higher dimension $M\ge
n$). The nonnegativity restriction on the entries of $\mathbf{A}$,
$\mathbf{b}$, $\mathbf{c}$ reflects physical constraints in
applications. Such positive systems appear, for example, in
modeling of bio-systems, chemical reaction systems, and
socio-economic systems, as described in detail in
\cite{farri,kacz,luen}. A thorough overview of the positive
realization problem and related results has recently been given in
\cite{tutorial}, while for a direct application in filter-design
we refer the reader to \cite{fibr}.

The {\sl existence problem} is to decide for a given transfer
function whether any positive realization $\mathbf{A}$,
$\mathbf{b}$, $\mathbf{c}$ of any dimension $M$ exists. It is 
known that the constraint of positivity may force the dimension $M$ to be strictly larger than $n$, see \cite{elozmenytolarge}, \cite{large}, \cite{mn2} for
different reasons why this phenomenon may occur.
The {\sl minimality problem} is to find the lowest possible value of $M$. These
problems have been given considerable attention over the past
decade. The existence problem was completely solved in \cite{negy}
and \cite{far}, cf., \cite{kad,mn1,fn}, while a few
particular cases of the minimality problem were settled in
\cite{fb2,halm,posmulti,astolfi,kinai,picci}.

The state of the art of the theory is therefore rather two-sided.
On one hand, there exists a {\sl general and constructive
solution} \cite{negy,far} to the existence problem which, however,
is {\sl inefficient} in the sense that it yields very large
dimensions, even in trivial cases. On the other hand, the
minimality problem is solved only for {\sl particular classes} of
transfer functions, and a general solution seems to be out of reach
for current methods. 


\smallskip

In this note we first observe (Section \ref{sec2}) that an {\sl
appropriate combination of known results leads to a constructive,
efficient, general algorithm} to solve the existence problem in
close-to-minimal dimensions. We observe that a repeated
application of a lemma of Hadjicostis \cite{had} leads to the 
{\sl positive decomposition problem} which, in turn, may be
treated by methods of \cite{dij,compmulti}. In Section \ref{sec3} we
give two illustrative examples. In the first we compare the
arising dimension to that of the earlier general algorithm of
\cite{negy}. In the second we complement the results of
\cite{large} by determining the minimal value of $M$ for a class
of transfer functions. Finally, in Section \ref{sec4} we provide a
{\it new lower-bound} on $M$, improving a result of \cite{mn2}.
This latter contribution is independent of earlier results.

\section{The algorithm}\label{sec2}


It is known that a necessary condition for the existence of positive realizations is that
one of the dominant poles (i.e. the poles
with maximal modulus) of $H(z)$ be nonnegative real, and
there is no loss of generality in assuming that it is located at
$\l_0=1$, see, e.g. \cite{negy}. The transfer function $H(z)$ is
called {\it primitive} if $\l_0$ is a unique dominant pole. It is also
known, see \cite{far}, that by the method of down-sampling the
case of non-primitive transfer functions can be traced back to
primitive ones. Therefore it is customary to assume that $H(z)$
is a {\it primitive transfer function with dominant pole at
$\l_0=1$}. We shall also assume, for technical simplicity, that
$\l_0=1$ is a {\it simple pole} (this makes the calculations less
involved; we note that the case of a multiple dominant pole can be
reduced to the simple pole case as in \cite[Step 4]{mn1}). Without loss of generality we may assume that the residue at $\l_0=1$ is 1 (see e.g. \cite{negy}).

With these normalizing assumptions, the transfer function
$H(z)$ takes the form
\begin{eqnarray}
H(z)&=&\frac{1}{z-1}+G(z)\nonumber\\
&=&\frac{1}{z-1}+\sum_{j=1}^r
\sum_{i=1}^{n_j}\frac{c_j^{(i)}}{(z-\l_j)^i},\label{alap}
\end{eqnarray}
where the poles $\l_j$ of $G(z)$ are of modulus strictly less than
1, i.e. $G(z)$ is asymptotically stable (note that 
$\l_j$s and $c_j^{(i)}$s are possibly complex).

\smallskip

In the series expansion $\displaystyle H(z)=\sum_{k=1}^\infty
t_kz^{-k}$ the coefficients $t_k$ are called the \emph{impulse
response} of $H(z)$.

\smallskip

If $H(z)={\bf{c}}^T(z{\bf{I}}-{\bf{A}})^{-1}{\bf{b}}$
then $t_k={\bf{c}}^T{\bf{A}}^{k-1}{\bf{b}}$ for all $k\geq0$. In
particular, $t_k$s must be non-negative for $H(z)$ to have a
positive realization. We now give the main ingredients upon which
the algorithm is based. The first is the following simple but
powerful result of Hadjicostis (see \cite[Theorem 5]{had}).

\smallskip

\begin{lemma}\label{lem1} (Hadjicostis)\ \\
{\sl Let $\displaystyle H(z)=\sum_{j=1}^\infty t_jz^{-j}$ be a rational transfer
function with \emph{non-negative} impulse response $t_1, t_2,
\dots$. For $m\ge 1$ let $H_m(z)$ denote the transfer
function corresponding to the shifted sequence $t_m, t_{m+1},
\dots$, i.e. $\displaystyle H_m(z)=\sum_{j=1}^\infty t_{m+j-1}z^{-j}$.
Assume that $H_m(z)$ admits a positive realization of some
dimension $k$. Then $H(z)$ admits a positive realization of
dimension $k+m-1$.}
\end{lemma}

\smallskip

We apply Lemma \ref{lem1} to $H(z)$ as given in (\ref{alap}). Note
that $H_1(z)=H(z)$ by definition, and for each $m\ge 2$ we have
$H_m(z)=zH_{m-1}(z)-t_{m-1}$. Hence, for each $m\ge 1$, 
$H_m(z)=\frac{1}{z-1}+\sum_{j=1}^r
\sum_{i=1}^{n_j}\frac{c_{j,m}^{(i)}}{(z-\l_j)^i}$. {\it The
leading coefficient remains 1}, while all other coefficients
$c_{j,m}^{(i)}\to 0$ exponentially as $m\to \infty$ (due to the
asymptotic stability of $G(z)$). That is, the leading coefficient
becomes {\it large compared to other coefficients}, and this
is exactly the familiar situation of the {\it positive
decomposition problem}, which we now turn to.

The task in the positive decomposition problem is to decompose an
arbitrary transfer function $G(z)$ as the difference
$G(z)=T_1(z)-T_2(z)$, with $T_1(z)$ and $T_2(z)$ both admitting
positive realizations (see \cite{dij,halm,posmulti,compmulti}). By
rescaling, one may assume that $G(z)$ is asymptotically stable, and
then the usual approach is to take a one-dimensional positive
system $T_2(z)=\frac{R}{z-w}$, where $0<w<1$ is larger than the
modulus of any pole of $G(z)$, and $R$ is a sufficiently large
positive number. Then $T_1(z)=G(z)+T_2(z)$ can be shown to admit a
positive realization which, in some cases, turns out to be also
minimal \cite{dij,halm,posmulti,compmulti}. For our purposes, the
essence of these results can be summarized as follows: {\it for
any primitive transfer function, as long as the partial fraction
coefficient of the dominant pole is significantly larger than all other
coefficients (as in $H_m(z)$ and $T_1(z)$ above) there exist
efficient methods to construct positive realizations.} We shall not list all relevant results of \cite{dij,halm,posmulti,compmulti} concerning the positive
decomposition problem; instead, we give as an example Theorem 8 of
\cite{dij}, which handles all transfer functions with simple poles.


\medskip

\begin{theorem}\label{simplepole} (Benvenuti, Farina \& Anderson)\ \\
{\sl Let $\displaystyle H(z)=\frac{1}{z-1}+G(z)=\frac{1}{z-1}+\sum_{j=1}^{n-1}
\frac{c_j}{z-\l_j}$, where $G(z)$ is a strictly proper
asymptotically stable rational transfer function of order $n$,
with {\it simple} poles. Let $\mathcal{P}_j$ $(j\ge 3)$ denote
the interior of the regular polygon in the complex plane with $j$ edges
located at the $j$-th roots of unity.
$\mathcal{P}_j$ can formally be defined in
polar coordinates as in \cite{dij}:
\begin{eqnarray*}
\mathcal{P}_j:=\Bigl\{(\rho ,\theta )\,:
 \ \rho
\cos\Bigl(\frac{(2k+1)\pi}{j}-\theta\Bigr)<
\cos\frac{\pi}{j},\ \quad\\
\mathrm{for}\ k=0,1,\dots ,j-1\Bigr\}.
\end{eqnarray*}
Let $N_1$ be the number of non-negative real poles with positive residue in $G(z)$
and let $N_2$ denote the number of other real poles 
in $G(z)$. Let $N_3$ denote the number of pairs
of complex conjugate poles of $G(z)$ belonging to the region
$\mathcal{P}_3$, and let $N_j$ $(j\ge 4)$ denote the number of pairs
of complex conjugate poles of $G(z)$ belonging to the region
$$
\displaystyle\mathcal{P}_j\setminus \bigcup_{m=3}^{j-1}\mathcal{P}_m .
$$
If all $c_j$s are sufficiently small then $H(z)$ admits a positive realization of dimension
$\displaystyle N=(n-1)+N_2+\sum_{j\ge 3}(j-2)N_j=\sum_{j\ge 1}jN_j$.}
\end{theorem}

\begin{center}
\begin{tabular}{ccccc}
\begin{xy}
<10mm,0mm>*{} ; <-5mm,8mm>*{} **\dir{-},
<-5mm,8mm>*{}; <-5mm,-8mm>*{} **\dir{-},
<-5mm,-8mm>*{}; <10mm,0mm>*{} **\dir{-},
\end{xy}&\ \ \ \ \ &
\begin{xy}
<10mm,0mm>*{} ; <-5mm,8mm>*{} **\dir{.},
<-5mm,8mm>*{}; <-5mm,-8mm>*{} **\dir{.},
<-5mm,-8mm>*{}; <10mm,0mm>*{} **\dir{.},
<10mm,0mm>*{} ; <0mm,10mm>*{} **\dir{-},
<0mm,10mm>*{} ; <-10mm,0mm>*{} **\dir{-},
<-10mm,0mm>*{} ; <0mm,-10mm>*{} **\dir{-},
<0mm,-10mm>*{} ; <10mm,0mm>*{} **\dir{-},
\end{xy}&\ \ \ \ \ &
\begin{xy}
<10mm,0mm>*{} ; <-5mm,8mm>*{} **\dir{.},
<-5mm,8mm>*{}; <-5mm,-8mm>*{} **\dir{.},
<-5mm,-8mm>*{}; <10mm,0mm>*{} **\dir{.},
<10mm,0mm>*{} ; <0mm,10mm>*{} **\dir{.},
<0mm,10mm>*{} ; <-10mm,0mm>*{} **\dir{.},
<-10mm,0mm>*{} ; <0mm,-10mm>*{} **\dir{.},
<0mm,-10mm>*{} ; <10mm,0mm>*{} **\dir{.},
<10mm,0mm>*{} ; <3mm,9mm>*{} **\dir{-},
<3mm,9mm>*{} ; <-8mm,6mm>*{} **\dir{-},
<-8mm,6mm>*{} ; <-8mm,-6mm>*{} **\dir{-},
<-8mm,-6mm>*{}; <3mm,-9mm>*{} **\dir{-},
<3mm,-9mm>*{};<10mm,0mm>*{} **\dir{-},
\end{xy}
\\
$\mathcal{P}_3$&&--- $\mathcal{P}_4$ &&
--- $\mathcal{P}_5$\\
&&$\cdots$ $\mathcal{P}_3$&& $\cdots$ $\mathcal{P}_3$ and $\mathcal{P}_4$
\end{tabular}
\smallskip

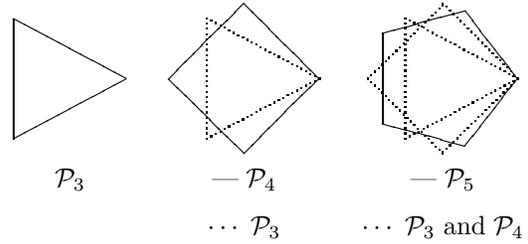
{\sc Figure 1.} {The sets $\mathcal{P}_j$.}
\end{center}

\medskip

\begin{remark}\label{r1}\ \\
The dimension $N$ appearing in the theorem is not necessarily
minimal but it is a good {\it a priori} upper bound on the
order of the realization. Further, by carefully analysing the proof in \cite{dij}, the condition on the residues may be given explicitely to be $\sum|c_j|\leq 2^{-5/2}$ (where the sum runs only over the $j$'s for which $\lambda_j$ is counted in $N_k$, $k\geq2$). In such case it can be computed that $H(z)$ admits a positive realization of dimension
$\displaystyle N=\sum_{j\ge 1}jN_j$.
For the reader's convenience, we have indicated in the appendix 

The proof in \cite{dij} is constructive and it gives a so-called \emph{cone-generated} positive realization (see Section 4).
\end{remark}
\medskip

This theorem was later improved and
generalized in various forms, \cite[Corollary 2]{halm},
\cite[Corollary 2]{posmulti}, \cite[Theorem 1, Theorem
2]{compmulti}. These papers also provide a number of examples
where {\it minimality} of the arising dimension $N$ can be
claimed. Finally, a synthesis of all these results, \cite[Theorem4]{compmulti}, covers the case of $H(z)=\frac{1}{z-1}+G(z)$ for
{\it any} asymptotically stable rational transfer function $G(z)$.

\medskip

\begin{theorem}\label{general}(Matolcsi, Nagy \& Szilv\'asi)\ \\
If $G(z)$ is any asymptotically stable rational transfer function
with poles $\l_1, \dots \l_r$ of order $m_1, \dots m_r$, and if all the
partial fraction coefficients of $G$ are sufficiently small, then the function
$H(z)=\frac{1}{z-1}+G(z)$ admits a positive realization, the
dimension of which is given explicitely as a function of $\l_1,
\dots \l_r, m_1, \dots m_r$.
\end{theorem}

\medskip

We remark that in the case of simple poles the value of $N$ in
Theorem \ref{simplepole} is better than in Theorem \ref{general}.
All the results above are constructive (see \cite{dij,halm,posmulti,compmulti}).

Theorems \ref{simplepole} and \ref{general} yield that the
following is a {\it general} algorithm which terminates
in a finite number of steps for any given transfer function $H(z)$.

\medskip

{\bf ALGORITHM:}

\medskip

{\it Assume $H(z)$ is given as in (\ref{alap}).} 


\noindent $m:=1$
\vskip0.1cm
\noindent WHILE $t_m \ge 0$ DO
\vskip0.1cm
\noindent \hskip0.5cm IF assumptions of Theorem \ref{simplepole} = TRUE 
\vskip0.1cm
\noindent \hskip0.5cm THEN APPLY Theorem \ref{simplepole}
\vskip0.1cm
\noindent \hskip1cm APPLY Lemma \ref{lem1}
\vskip0.1cm
\noindent \hskip0.5cm ELSE 
\vskip0.1cm
\noindent \hskip0.5cm IF assumptions of Theorem \ref{general} = TRUE 
\vskip0.1cm
\noindent \hskip0.5cm THEN APPLY Theorem \ref{general}
\vskip0.1cm
\noindent \hskip1cm APPLY Lemma \ref{lem1}
\vskip0.1cm
\noindent \hskip0.5cm ELSE $m:= m+1$
\vskip0.1cm
\noindent ELSE there is no positive realization of $H(z)$

\medskip

\begin{remark}\label{rem1}\ \\
At each step, we may apply to $H_m(z)$
{\it other known constructions} from the literature different than
Theorems \ref{simplepole} or \ref{general}. We included only these
two theorems to keep the algorithm transparent and because they
guarantee that the algorithm terminates in a finite number of steps. Other
important partial results which can be incorporated into the
algorithm are given in \cite{fb2,astolfi,picci}.
\end{remark}

\medskip

\begin{remark}\label{re2}\ \\
We acknowledge that this algorithm is merely
an observation that some earlier results in the literature can be
combined together. Nevertheless, we find it an important observation
as it provides a {\it completely general} algorithm. Previously
such a general algorithm has only been given in \cite{negy}.
There are good heuristic arguments to believe that
the algorithm above is efficient in terms of producing small
dimensions, and better than the existing algorithm of
\cite{negy}. First, the partial fraction coefficients decay
exponentially, so that only a few iterations are needed before
Theorem \ref{simplepole} or \ref{general} become applicable, and
these theorems already provide minimal or close-to-minimal
dimensions (see Section \ref{sec3} for a numerical examples).
Second, the method of \cite{negy} involves the time development of
an $n$-dimensional "cube" around the vector $(1, 1, \dots, 1)$
and, as such, can only produce dimensions larger than $2^n$
(usually significantly  larger than that).  This fact, however, by
no means diminishes the theoretical significance of the results of
\cite{negy} which provided the first general solution to the {\it
existence problem} for primitive transfer functions.
\end{remark}

\section{Examples}\label{sec3}

In this section we give two examples. In the first we compare the
arising dimension of realization with that of the algorithm of
\cite{negy}. In the second we complement a result of \cite{large}
and determine the minimal dimension of positive realizations for a
class of transfer functions.

First, we note that in the case where there are only simple poles,
the number of iterations needed may be evaluated as follows:
write $H(z)=\sum_{j=0}^n\frac{c_j}{z-\lambda_j}$ with $\lambda_0=1$ and
$|\lambda_j|<1$. A simple computation shows that
$H_m(z)=\sum_{j=0}^n\frac{c_j\lambda_j^{m-1}}{z-\lambda_j}$. It follows that,
if $m\geq \left|\frac{\log 2^{5/2}n\max|c_i|}{\log\max|\lambda_i|}\right|$, then
$\sum |c_j\lambda_j^{m-1}|\leq 2^{-5/2}$ and Theorem \ref{simplepole} applies (cf. Remark \ref{r1}).
Moreover, it is enough to consider those poles that are not non-negative with non-negative residues.

\medskip

\noindent{\bf Example 1.} Let
\begin{equation*}H(z)=\frac{1}{z-1}+t(z)=\frac{1}{z-1}+
\end{equation*}
\begin{equation*}
{\frac {0.3331328522\,{z}^{2} +0.1984152016\,z+0.1253986950}{{z}^{3}- 0.69055619\,{z}^{2}+ 0.80189061\,z-
0.38920832}},
\end{equation*}
where $t(z)$ is a low-pass digital Chebyshev filter of order 3.
The partial fraction decomposition is
\begin{equation*}
H(z)=\frac{1}{z-1}+\sum_{i=1}^3\frac{c_i}{z-\lambda_i}
\end{equation*}
with
\begin{eqnarray*}
\lambda_1&=&0.07522998673 - 0.8455579204\,i,\quad\lambda_2=\overline{\lambda_1}\\
c_1&=&- 0.01050864690+0.1411896961\,i,\quad c_2=\overline{c_1}\\
\lambda_3&=&0.5400962165\quad
c_3=0.3541501460.
\end{eqnarray*}
Let us now apply the algorithm. The above computation shows that after $m=5$ iterations, we may
apply Theorem \ref{simplepole} to $H_5$. As $N_2=0$ and as 
$\lambda_1,\lambda_2\in\mathcal{P}_4\setminus\mathcal{P}_3$, $N_4=1$, $N_j=0$ for $j\not=4$
we obtain a positive realization of $H_5(z)$ of dimension $N=5$ (cf. Remark \ref{r1}). A repeated application of
Lemma \ref{lem1} gives a positive realization of
$H(z)$ of dimension $9$. A more careful examination of the
proof of Theorem \ref{simplepole} would show that
$H_3$ already satisfies the requirements so that a $7$-dimensional
positive realization exists. We spare the reader the numerical values.
%

While we cannot claim that the arising dimension is minimal let us
compare it with that of the algorithm described previously in
\cite{negy}. In short, a specific minimal realization of $H(z)$ is
considered and the time-evolution of a small 4-dimensional cube
around the vector $[1,1,1,1]^T$ needs to be checked to provide a
system-invariant polyhedral cone, which then leads to a positive
realization of $H(z)$. The number of edges of the cone equals the
dimension of the positive realization. If one carries out this
procedure word-by-word for $H(z)$ above the arising dimension is
$48$.

\medskip

\noindent{\bf Example 2.}
Consider the family of transfer functions
\begin{equation}\label{hn}
H^{N}(z)=\frac{1}{z-1}-\frac{4\cdot
(5/2)^{N-2}}{z-0.4}+\frac{3\cdot 5^{N-2}}{z-0.2}
\end{equation}
as in \cite[Example 4]{large}. It is proved in \cite{large} that
for any $N\ge 4$ the minimal dimension of positive realizations of
$H^N(z)$ is {\it at least} $N$. Here we prove that an
$N$-dimensional minimal positive realization of $H^N(z)$ does
indeed exist for every $N\ge 4$.

For $N=4$ the following 4-dimensional positive realization of
$H^4(z)$ is given in \cite[Example 3]{large} and it is shown to
be minimal.
$$
\mathbf{b}=\begin{pmatrix}0&0&0&1\end{pmatrix}^T,\quad
\mathbf{c}=\begin{pmatrix}6&0&0&51\end{pmatrix}^T
$$
and
\begin{eqnarray}\label{4dim}
&\mathbf{A}=\left(\begin{array}{cccc}
 0 & 0 & 0 & 1 \\
 1 & \frac{63+4\sqrt{26}}{85} & 0 & 0 \\
 0 & \frac{22-4\sqrt{26}}{85} & \frac{63-4\sqrt{26}}{85} & 0 \\
 0 & 0 & \frac{22+4\sqrt{26}}{85} & 0
\end{array}\right), \quad
&
\end{eqnarray}

Consider now $H^N(z)$, $N\ge 4$. It is not difficult to see that
with the notation of Lemma \ref{lem1} we have $H_1(z)=H^N(z)$,
$H_2(z)=H^{N-1}(z)$, $\dots$, $H_m(z)=H^{N+1-m}(z)$. Let us stop
at $m=N-3$, i.e. at $H_m(z)=H^4(z)$ and make use of the
realization (\ref{4dim}) of $H^4(z)$. Then, the application of
Lemma \ref{lem1} produces a positive realization of $H^N(z)$ of
order $4+(m-1)=N$. Note that for producing this minimal positive
realization we use (\ref{4dim}) instead of Theorem \ref{simplepole} (we
have, in fact, followed the suggestion of Remark \ref{rem1}).

We have concluded that the minimal dimension of positive
realization of $H^N(z)$ is $N$. Let us now see what the
word-by-word application of the algorithm of Section \ref{sec3}
gives.

The algorithm terminates when Theorem \ref{simplepole} becomes
applicable, i.e. when $H_m(z)=H^0(z)$, that is $m=N+1$. Then a
3-dimensional positive realization of $H^0(z)$ is constructed, and
the application of Lemma \ref{lem1} produces a positive
realization of $H^N(z)$ of order $3+(m-1)=N+3$.

This example is reassuring in that the application of the
algorithm of Section \ref{sec3} produces positive realizations
of close-to-minimal order.

\section{Improved lower-bounds}\label{sec4}

We saw in Section \ref{sec3} that the minimal order of positive
realizations of $H^N(z)$ is $N$. Also, it is easy to calculate
(see \cite{large}) that the impulse response sequence of $H^N(z)$
contains zeros, namely $t_{N-1}=t_N=0$. A general lower-bound
presented in \cite{mn2} gives that in such case the order $M$ of
any positive realization satisfies $\frac{M(M+1)}{2}-1+M^2\ge
N$, i.e. $M$ is at least $\approx \sqrt{\frac{2N}{3}}$. In view of
the actual minimal value $M=N$ a lower-bound of the order of
magnitude $N$ is welcome (instead of the order of
magnitude $\sqrt{N}$). In this section we present such an
improvement (but we note that while the lower-bound of
\cite{mn2} is valid in general, our improvement is restricted
to transfer functions with positive real poles, as is the case of
the example of the previous section).

Throughout this section we assume that $H(z)$ is a given primitive
transfer function of McMillan degree $n$ with {\it positive real
poles}, and there exists a positive integer $k_0$, such that for
the impulse response sequence of $H(z)$ we have $t_{k_0}=0$ and
$t_{k}>0$ for all $k>k_0$. This means that $H(z)$ is of the form
\begin{equation}\label{real}
\frac{1}{z-1}+\sum_{j=1}^r
\sum_{i=1}^{n_j}\frac{c_j^{(i)}}{(z-\l_j)^i}
\end{equation}
 where
$ c_j^{(i)}\in {\R}$, $0<\l_j<1$, and $\sum_{j=1}^{r}n_j=n-1$.

Let the triple
$(\mathbf{h},\mathbf{F},\mathbf{g})$ denote an arbitrary
minimal ($n$-dimensional) realization of $H(z)$ (for canonical
minimal realizations see e.g. \cite{farri}). Assume that there
exists a matrix $\mathbf{P}$ of size $n\times M$ such that for some
triple $({\bf{c}},{\bf{A}},{\bf{b}})$ with nonnegative entries:
\begin{eqnarray}\label{kup}
\mathbf{FP}=\mathbf{PA}, \ \  \mathbf{Pb}=\mathbf{g}, \ \
\mathbf{c}^T=\mathbf{h}^T\mathbf{P}.
\end{eqnarray}
There is a well-known geometrical interpretation of these
equalities.  Namely, the columns of matrix $\mathbf{P}$ represent
the edges of a finitely generated cone $\mathcal{P}$ in $\R^n$,
such that $\mathcal{P}$ is $\mathbf{F}$-invariant, and
$\mathcal{P}$ lies between the reachability cone and the
observability cone corresponding to the triple
$(\mathbf{h},\mathbf{F},\mathbf{g})$. It is known that the
triple $(\mathbf{c},\mathbf{A},\mathbf{b})$ provides a
positive realization of $H(z)$.

\medskip

\begin{definition}\ \\
{\sl A triple $(\mathbf{c},\mathbf{A},\mathbf{b})$
which arises in such a manner is called a
{\it{cone-generated}} realization of $H(z)$.}
\end{definition}

\medskip

It is a basic result in the theory of positive realizations that a
transfer function $H(z)$ admits positive realizations if and only
if it admits cone-generated realizations (see \cite{eredeti}). 
Here we present a lower bound on the order
of cone-generated realizations of $H(z)$. For this we shall need the
following auxiliary result.

\medskip

\begin{lemma}\label{lem2}\ \\
{\sl Let $f:\R\to\R$ be defined by
$$
f(x):=\sum_{j=1}^r p_{(j)}(x) \lambda_{j}^{x}
$$
where $\lambda_1>\lambda_2>\dots >\lambda_r>0$ and $p_{(j)}$ denotes a polynomial
(with real coefficients) of degree $n_j$. Then $f$ has at most
$R\equiv R_f:=\displaystyle\sum_{j=1}^{r}(n_j+1)-1$ pairwise distinct real
roots.}
\end{lemma}

\bigskip


Lemma \ref{lem2} is proved by induction on $R$.

We are now ready to give an improvement of the lower-bound of
\cite{mn2}.

\medskip

\begin{theorem}\label{theo2}\ \\
{\sl Assume that $H(z)$ is a transfer function of McMillan degree $n$,
with positive real poles, given as in (\ref{real}) above. Assume
also that there exists a positive integer $k_0$, such that for the
impulse response sequence of $H(z)$ we have $t_{k_0}=0$ and
$t_{k}>0$ for all $k>k_0$. Then the dimension $M$ of any
cone-generated positive realization of $H(z)$ satisfies $M\ge
\frac{k_0}{n-1}$.}
\end{theorem}

\bigskip

{\it Proof.}\quad
Let the triple $(\mathbf{h},\mathbf{F},\mathbf{g})$ denote a 
minimal ($n$-dimensional) realization of $H(z)$. Consider any cone-generated positive realization $({\bf{c}},{\bf{A}},{\bf{b}})$ of $H(z)$ arising from a matrix
$\mathbf{P}$ of size $n\times M$, as in \eqref{kup}.  Let ${\bf{e}}_i$ denote an
arbitrary column of the matrix $\mathbf{P}$, and consider the sequence
$g^{(i)}_k:=\mathbf{h}^T\mathbf{F}^{k-1}{\bf{e}}_i\ge 0$. Let
$\mathbf{P}_1:=[p_{i,j}]$ be the nonnegative matrix of size $M\times
\infty$ defined by $p_{i,j}:=g^{(i)}_j$ for $1\le i\le M$ and
$1\le j$.

Let $\mathbf{K}:=[k_{i,j}]$ denote
the infinite Hankel matrix composed of the impulse response
sequence of $H(z)$, i.e. $k_{i,j}:=t_{i+j-1}$.  By assumptions imposed  on $\mathbf{P}$ there exists
a matrix $\mathbf{Q}=[q_{i,j}]$ of size $\infty \times M$, with
nonnegative entries, such that ${\bf{QP}}_1=\mathbf{K}$. This
is true because the $k$th row of $\mathbf{K}$ is given by
$k_{k,j}=t_{k+j-1}=\mathbf{h}^T\mathbf{F}^{j-1}(\mathbf{F}^{k-1}\mathbf{g})$,
and the vector $\mathbf{F}^{k-1}\mathbf{g}$ lies inside the cone
$\mathcal{P}$ by assumption. Thus, it may be decomposed as a
linear combination of the edges ${\bf{e}}_i$ of $\mathcal{P}$
with nonnegative coefficients, and [one choice of] these
coefficients form the $k$th row of the matrix  $\mathbf{Q}$.

Since $(\mathbf{h},\mathbf{F},\mathbf{g})$ is a minimal realization, for
an arbitrary column ${\bf{e}}_i$ of the matrix $\mathbf{P}$, the
transfer function corresponding to the impulse response sequence
$g^{(i)}_k=\mathbf{h}^T\mathbf{F}^{k-1}{\bf{e}}_i$ is of the form
$$
H^{({\bf{e}}_i)}(z)=\frac{C_i}{z-1}+\sum_{j=1}^r
\sum_{s=1}^{n_j}\frac{d_{j,i}^{(s)}}{(z-\l_j)^s}.
$$
(Note that some coefficients $C_i$ and $d_{j,i}^{(s)}$ may be 0.)  The {\it
column} ${\bf{e}}_i$ of $\mathbf{P}$ is called {\it{dominant}}
if $C_i\ne 0$ in $H^{({\bf{e}}_i)}(z)$.
Delete the non-dominant rows from the matrix $\mathbf{P}_1$ and
the corresponding columns from the matrix $\mathbf{Q}$. The
remaining matrices (of sizes $M_1\times\infty$ and $\infty\times
M_1$ for some $M_1\le M$) are denoted by
$\mathbf{P}_{1}^{(dom)}$ and $\mathbf{Q}^{(dom)}$. We see that
$\mathbf{Q}^{(dom)}\mathbf{P}_{1}^{(dom)}\le \mathbf{K}$
entrywise. Recall that $t_{k_0}=0$, by assumption on
the impulse response of $H(z)$. This implies that for some
dominant index $i$ ($1\le i\le M_1$),
$g^{(i)}_{k_0}=0$. Otherwise, $k_{1,k_0}=t_{k_0}$ would be strictly
positive in the first row of $\mathbf{K}$. Considering the
second row of $\mathbf{K}$ we see that $k_{2,k_0-1}=t_{k_0}=0$,
hence $g^{(i)}_{k_0-1}=0$ for some dominant index
$i$. By the same argument, for every $1\le j\le k_0$,
$g^{(i)}_{j}=0$ for some dominant index $i$. In other
words, each of the first $k_0$ columns of the matrix
$\mathbf{P}_{1}^{(dom)}$ contains a zero, and hence there are
at least $k_0$ zero entries in $\mathbf{P}_{1}^{(dom)}$.

On the other hand, $g^{(i)}_k$ is the impulse response of
$H^{(e_i)}(z)=\frac{C_i}{z-1}+\sum_{j=1}^r
\sum_{s=1}^{n_j}\frac{d_{j,i}^{(s)}}{(z-\l_j)^s}$. Thus,
$g^{(i)}_k=C_i+\sum_{j=1}^{r-1} p_{(j)}(k) \l_{j}^{k}$, where
$p_{(j)}$ are polynomials of degree not exceeding $n_j-1$, for
$1\le j\le r$. Hence, Lemma \ref{lem2} implies that there are at
most $R=(1+\sum_{j=1}^{r}n_j)-1=n-1$ zeros in each row of
$\mathbf{P}_{1}^{(dom)}$. This means that the number of zeros in the
matrix $\mathbf{P}_{1}^{(dom)}$ is at most $M_1(n-1)$. Therefore,
$$
k_0\le \#\mbox{ of zeros in }\mathbf{P}_{1}^{(dom)}\le M_1(n-1)\le M(n-1),
$$
and, hence, $M\ge \frac{k_0}{n-1}$.
\hfill$\Box$

\bigskip

\begin{remark}\ \\
If we apply Theorem \ref{theo2} to functions
$H^N(z)$ of Section \ref{sec3} we obtain $M\ge N/2$. This is still
far from the actual minimal value $N$. However, if there are
more than 3 poles present in $H(z)$ then the geometric arguments
of \cite{large} are difficult to generalize, while Theorem
\ref{theo2} still applies. 
\end{remark}

\bigskip

\begin{remark}\ \\
As mentioned in the ``Open Problems and New Directions" section of \cite{tutorial} it is
desirable to have tight upper and lower bounds on the minimal
order of a positive realization {\it in general}. Note, however,
that the results of \cite{large,mn2} and Theorem \ref{theo2} above
are all based on the assumption that the impulse response sequence
of $H(z)$ contains at least one 0. The only other
lower-bound known to us is that of \cite{had} which, however, does
not give any non-trivial estimates for transfer functions with
nonnegative poles.

What can be said if the impulse response does not contain zeros?
Unfortunately, we do not have a general approach to this case. As
a first step in this direction we examined the modified family
$H^{N,\e}(z)=\frac{1}{z-1}-\frac{4\cdot
(5/2)^{N-2}}{z-0.4}+\frac{3\cdot 5^{N-2}+\e}{z-0.2}$ for small
values of $\e$. Note that the impulse response sequence no longer
contains zeros. Since the dimension of the system is 3, we can use
elementary (but tedious) geometric arguments to conclude that for
small enough $\e$ the minimal order $M$ of positive realizations
of $H^{N,\e}(z)$ still satisfies $M\ge N/2$. It is not clear,
however, how to generalize these arguments to transfer functions
of higher degree (as in Theorem \ref{theo2}) where the geometric
intuition is missing. Therefore, finding tight lower-bounds in the
general case remains an open problem.
\end{remark}

\section{Conclusion}

We have observed that recent results in positive system theory can
be put together to produce an {\it efficient, general} algorithm
to the positive realization problem of transfer functions. We have
given two examples to illustrate the algorithm. In the first we
compared the arising dimension of realization with that of an
earlier general algorithm of \cite{negy}. In the second we
examined a family of transfer functions given in \cite{large}, and
determined the minimal order of positive realizations. With
respect to the minimality problem we have proved a new lower-bound
on the order of positive realizations of transfer functions with
positive real poles, improving an earlier general result of
\cite{mn2}.

\appendix

\section{Precisions on the paper \cite{dij}}

In this section, we will show how to obtain a quantitative bound on $R$ in
Theorem \ref{simplepole} above. This result is not really new as it can be
obtained directly from the proof of that theorem in \cite{dij} and some
simple observations. This appendix is thus only included here for the reader's convenience.

\begin{proposition}{\cite[Proposition 7]{dij}}
Let $H(z)=\frac{\eta e^{i\vartheta}}{z-\rho e^{i\theta}}+\frac{\eta e^{-i\vartheta}}{z-\rho e^{-i\theta}}$
and assume that $\rho e^{i\theta}\in\mathcal{P}_m$. Then for 
$\displaystyle R\geq \frac{2\eta}{\cos\frac{\pi}{m}}$
there exists $A_+\in\R_+^{m\times m}$, $b_+,c_+\in\R_+^m$ such that
$$
H_1(z)=H(z)+\frac{R}{z-1}=c_+^T(zI-A_+)^{-1}b_+
$$
\end{proposition}

\noindent{\it Proof.}
Let us first consider the Jordan realization of $H_1$\,: $H_1(z)=c^T(zI-A)^{-1}b$ with
$$
A=\begin{pmatrix}
\rho\cos\theta&-\rho\sin\theta&0\\
\rho\sin\theta&\rho\cos\theta&0\\
0&0&1 \end{pmatrix},
$$
$$
b=\begin{pmatrix}2\eta (\cos\vartheta+\sin\vartheta)\\ 2\eta (\cos\vartheta-\sin\vartheta)\\ R\end{pmatrix},
\quad  c=\begin{pmatrix}1\\ 1\\ 1\end{pmatrix}.
$$
Next, define $K$ by
$$
\begin{pmatrix}
1&\cos\left(\frac{2\pi}{m}\right)&\cos\left(\frac{2\pi}{m}2\right)&\cdots
&\cos\left(\frac{2\pi}{m}(m-1)\right)\\
0&\sin\left(\frac{2\pi}{m}\right)&\sin\left(\frac{2\pi}{m}2\right)&\cdots
&\sin\left(\frac{2\pi}{m}(m-1)\right)\\
1&1&1&\cdots&1\\
\end{pmatrix}
$$
and, for $\alpha>0$ define
$$
D(\alpha)=\begin{pmatrix} \alpha&0&0\\ 0&\alpha&0\\ 0&0&1\end{pmatrix}.
$$
It was proved in \cite{dij} that the positive cone $\mathcal{K}_\alpha$ generated by $D(\alpha)K$
is $A$-invariant. To prove the theorem, it is then enough to prove that we can choose $\alpha$
in such a way that
$$
\mathcal{K}_\alpha\subset\mathcal{O}:=\{\mathbf{x}\in\R^3\,: c^TA^k\mathbf{x}\geq 0,\quad k=0,1,\ldots\}
$$
and that if $R\geq \frac{2\eta}{\cos\frac{\pi}{m}}$ then $b\in\mathcal{K}_\alpha$.

For the first one, note that if 
$$
\mathbf{x}=\lambda(\alpha\cos\varphi,\alpha\sin\varphi,1)^T\in\R^3
$$
then
\begin{eqnarray*}
&&c^TA^k\mathbf{x}\\
&=&\lambda(1,1,1)\begin{pmatrix}
\rho^k\cos k\theta&-\rho^k\sin k\theta&0\\
\rho^k\sin k\theta&\rho^k\cos k\theta&0\\
0&0&1 \end{pmatrix}\begin{pmatrix}\alpha\cos\varphi\\ \alpha\sin\varphi\\ 1 \end{pmatrix}\\
&=&\lambda\bigl[\alpha\rho^k\bigl(\cos(\varphi+k\theta)+\sin(\varphi+k\theta)\bigr)+1\bigr].
\end{eqnarray*}
It follows that, for $\alpha=1/2$, $\lambda\geq0$, $\rho\leq 1$, we always have
$c^TA^k\mathbf{x}\geq 0$, in particular, $\mathcal{K}_{1/2}\subset\mathcal{O}$.

Finaly, note that the cone $\mathcal{K}_\alpha$ contains the cone
$\{\lambda(r,\alpha)\,:\ r<\cos\frac{\pi}{m},\lambda>0\}$, in particular
$b\in\mathcal{K}_\alpha$ if
\begin{equation}
\label{eq:djpp7}
\frac{2^{3/2}\eta}{R\alpha} <\cos\frac{\pi}{m}.
\end{equation}
By taking $\alpha=1/2$
we can see that this is the case as soon as $R\geq \frac{2^{5/2}\eta}{\cos\frac{\pi}{m}}$.
\hfill$\Box$

If the filter is given and that one seeks for more precise estimates,
one can slightly improve the result by taking $\alpha=1/2\rho$ and still a bit further if 
$\theta=r\pi$ for some $r$ (necessarily rational) for which one may have
$\cos(j\frac{2\pi}{m}+k\theta)+\sin(j\frac{2\pi}{m}+k\theta)\leq\kappa<1$ for
all integers $k\geq 0$ and $0\leq j\leq m$. In this case one may take $\alpha=\frac{1}{2\rho\kappa}$
and then check for the smallest $R$ for which
$\frac{2\eta}{R\alpha}(\cos\vartheta+\sin\vartheta,\cos\vartheta-\sin\vartheta)\in\mathcal{P}_m$.

In the opposit direction, note that $R=2^{7/2}\eta$ works for all $m\geq 3$.

Let us now show how the estimate on $R$ in Theorem \ref{simplepole} results from this.
We will take the notations from \cite{dij}. Let $q$ be the smallest integer so that
all complex poles of $H$ are in $\mathcal{P}_q$ and let us write
$$
H(z)=\sum_{j=1}^{N_1}\frac{c_j^{(1)}}{z-\lambda_j^{(1)}}+\sum_{i=2}^q\sum_{j=1}^{N_i}\left(
H_j^{(i)}(z)+\frac{R_j^{(i)}}{z-1}\right)
$$
where $c_j^{(1)}\geq0$ and $\lambda_j^{(1)}\geq0$, $H_j^{(2)}=\frac{c_j^{(2)}}{z-\lambda_j^{(2)}}$
corresponds to the other real poles and, for $i\geq 3$
$$
H_j^{(i)}=\frac{\eta_j^{(i)}e^{i\vartheta_j^{(i)}}}{z-\rho_j^{(i)} e^{i\theta_j^{(i)}}}
+\frac{\eta_j^{(i)}e^{-i\vartheta_j{(i)}}}{z-\rho_j^{(i)} e^{-i\theta_j^{(i)}}}
$$
where $\rho_j^{(i)}e^{-i\theta_j^{(i)}}\in\mathcal{P}_i\setminus\bigcup{k\leq i}\mathcal{P}_k$.

Now, each $\frac{c_j^{(1)}}{z-\lambda_j^{(1)}}$ has a one dimensional
one-dimensional positive realization
$A_j^{(1)},b_j^{(1)},c_j^{(1)}$. According to \cite[Proposition 5]{dij},
if $R_j^{(2)}\geq |c_j^{(2)}|$, $H_j^{(2)}+\frac{R_j^{(2)}}{z-1}$ has a
two-dimensional positive realization $A_j^{(2)},b_j^{(2)},c_j^{(2)}$
(the estimate is clear from the end of the proof of that proposition in \cite{dij} and the restriction
$\Lambda<1$ is irrelevant).
According to the improvement of \cite[Proposition 7]{dij} given above, if $R_j^{(i)}\geq 4\eta_j^{(2)}$,
then $H_j^{(i)}+\frac{R_j^{(i)}}{z-1}$ has a $i$-dimensional positive realization $A_j^{(i)},b_j^{(i)},c_j^{(i)}$. A positive realization dimension $N$ (defined in Theorem \ref{simplepole}) of $H$ is then given by 
\begin{eqnarray*}
b_+&=&\bigl(b_1^{(1)},\ldots,b_{N_1}^{(1)},b_1^{(2)},\ldots,b_{N_2}^{(2)},\ldots\\
&&\phantom{\bigl(b_1^{(1)},\ldots,b_{N_1}^{(1)},b_1^{(2)},\ldots}\ldots,
b_1^{(q)},\ldots,b_{N_q}^{(q)}\bigr)^T,\\
c_+&=&\bigl(c_1^{(1)},\ldots,c_{N_1}^{(1)},c_1^{(2)},\ldots,c_{N_2}^{(2)},\ldots\\
&&\phantom{\bigl(b_1^{(1)},\ldots,b_{N_1}^{(1)},b_1^{(2)},\ldots}\ldots,
c_1^{(q)},\ldots,c_{N_q}^{(q)}\bigr)^T
\end{eqnarray*}
and $A_+$ is given in block-diagonal notation by
\begin{eqnarray*}
A_+&=&\mbox{diag}\bigl(A_1^{(1)},\ldots,A_{N_1}^{(1)},A_1^{(2)},\ldots,A_{N_2}^{(2)},\ldots\\
&&\phantom{\bigl(b_1^{(1)},\ldots,b_{N_1}^{(1)},b_1^{(2)},\ldots}\ldots,
A_1^{(q)},\ldots,A_{N_q}^{(q)}\bigr)^T,
\end{eqnarray*}
A condition in the theorem is thus that $R:=\sum R_j^{(i)}$
satisfies {\it e.g.} $R\geq \sum_{j=1}^{N_2} |c_j^{(2)}|+2^{7/2}\sum_{i=3}^q\sum_{j=1}^{N_i}\eta_j^{(i)}$.
This bound may be improved slightly with the above remark.

Note also that this realization has dimension $N_1+2N_2+3N_3\cdots+qN_q$.

\bibliographystyle{IEEE}

\end{document}